\documentclass[11pt,reqno]{amsart}
\usepackage{amsmath, amssymb, amsthm}
\usepackage{url}
\usepackage[breaklinks]{hyperref}
\usepackage{array}
\renewcommand{\(}{\left\(}
\renewcommand{\)}{\right\)}
\renewcommand{\[}{\left\[}
\renewcommand{\]}{\right\]}

\numberwithin{equation}{section}
\theoremstyle{plain}
\newtheorem{theorem}{Theorem}[section]
\newtheorem{lemma}[theorem]{Lemma}
\newtheorem{remark}[]{Remark}

\newtheorem{conjecture}[theorem]{Conjecture}

\makeatletter
\def\proof{\@ifnextchar[{\@oproof}{\@nproof}}
\def\@oproof[#1][#2]{\trivlist\item[\hskip\labelsep\textit{#2 Proof of\
		#1.}~]\ignorespaces}
\def\@nproof{\trivlist\item[\hskip\labelsep\textit{Proof.}~]\ignorespaces}

\makeatother

\newenvironment{proof-alt}
{\vskip 0.15in \par\noindent{\it Proof of Proposition \ref{HSS}.}\hskip 0.5em\ignorespaces}
{\hfill $\Box$\par\medskip}

\begin{document}
	\title[Sign Patterns and Congruences]{Sign Patterns and Congruences of certain infinite products involving the Rogers-Ramanujan continued fraction}
	
	\author{Nayandeep Deka Baruah}
	\address{Nayandeep Deka Baruah, Department of Mathematical Sciences, Tezpur University, Napaam, Assam, 784028, India.}
	\email{nayan@tezu.ernet.in, nayandeeptezu@gmail.com}
	
	\author{Abhishek Sarma}
	\address{Abhishek Sarma, Department of Mathematical Sciences, Tezpur University, Napaam, Assam, 784028, India.}
	\email{abhitezu002@gmail.com}
	
	\thanks{2020 \textit{Mathematics Subject Classification.} 11P81, 11P83, 05A17.\\
		\textit{Keywords and phrases. Rogers-Ramanujan functions, congruences, periodicity.}}
	
\begin{abstract}
		We study the behavior of the signs of the coefficients of certain infinite products involving the Rogers-Ramanujan continued fraction. For example, if $$\sum_{n=0}^{\infty}A(n)q^{n}:= \dfrac{(q^2;q^5)_\infty^5(q^3;q^5)_\infty^5}{(q;q^5)_\infty^5(q^4;q^5)_\infty^5},$$then
		$A(5n+1)>0$, $A(5n+2)>0$, $A(5n+3)>0$, and $A(5n+4)<0$. 
		We also find a few congruences satisfied by some coefficients. For example, for all nonnegative integers $n$, $A(9n+4)\equiv 0 \pmod3$, $		A(16n+13)\equiv 0 \pmod4$, and $A(15n+r)\equiv0\pmod{15}$, where $r\in\{4, 8, 13, 14\}$.
\end{abstract}
	\maketitle
\section{Introduction}
The celebrated Rogers-Ramanujan continued fraction is defined by 
$$\mathcal{R}(q):=\dfrac{q^{1/5}}{1}_{+}\dfrac{q}{1}_{+}\dfrac{q^2}{1}_{+}\dfrac{q^3}{1}_{+~\cdots},\quad|q|<1.$$
Set $$R(q):=q^{-1/5}\mathcal{R}(q):=\dfrac{1}{1}_{+}\dfrac{q}{1}_{+}\dfrac{q^2}{1}_{+}\dfrac{q^3}{1}_{+~\cdots}.$$
The Rogers-Ramanujan identities are given by
\begin{align*}
	G(q)=\sum_{n=0}^{\infty}\frac{q^{n^2}}{(q;q)_n}=\frac{1}{(q;q^5)_{\infty}(q^4;q^5)_{\infty}}\\\intertext{and}
	H(q)=\sum_{n=0}^{\infty}\frac{q^{n^2+n}}{(q;q)_n}=\frac{1}{(q^2;q^5)_{\infty}(q^3;q^5)_{\infty}},
\end{align*}
where throughout the sequel, for complex numbers $a$ and $q$ with $|q|<1$, we use the customary notation
\begin{align*}
	(a;q)_0&:=1,~(a;q)_n:=\prod_{k=0}^{n-1}(1-aq^k) ~\textup{for}~ n\geq1,~\textup{and}~(a;q)_\infty:=\lim_{n\to \infty}(a;q)_n.
\end{align*}
Rogers \cite{rogers1894} and Ramanujan \cite{nb} (see \cite[Corollary, p. 30]{bcb3}) proved that  
\begin{align}\label{R-GH}R(q)&=\dfrac{H(q)}{G(q)}=\dfrac{(q;q^5)_\infty(q^4;q^5)_\infty}{(q^2;q^5)_\infty(q^3;q^5)_\infty}.
\end{align}

In 1978, Richmond and Szekeres \cite{richmond1978taylor} examined asymptotically the power series coefficients of a large class of infinite products including the product given in \eqref{R-GH} and its reciprocal. In particular, they \cite[Eq. (3.9)]{richmond1978taylor} proved that, if 
\begin{align*}
	\frac{1}{R(q)}=\dfrac{(q^2;q^5)_\infty(q^3;q^5)_\infty}{(q;q^5)_\infty(q^4;q^5)_\infty}=:\sum_{n=0}^{\infty} c(n)q^n,	
\end{align*}
then
$$c(n)=\dfrac{\sqrt{2}}{(5n)^{3/4}}\textup{exp}\left(\dfrac{4\pi}{25}\sqrt{5n}\right)\times\left\{\cos\left(\dfrac{2\pi}{5}\left(n-\dfrac{2}{5}\right)\right)+\mathcal{O}(n^{-1/2})\right\},$$
which implies that, for $n$ sufficiently large,
\begin{align}\label{thm:rich}
	c(5n)>0, ~c(5n+1)>0, ~c(5n+2)<0,~c(5n+3)<0,~\text{and}~ c(5n+4)<0.   
\end{align}
They also gave a similar result for the power series
coefficients of $R(q)$ from which it follows that, for $n$ sufficiently large,\begin{align}\label{thm:rich-2}
	d(5n)>0,~d(5n+1)<0,~ d(5n+2)>0,~d(5n+3)<0,~ \text{and}~d(5n+4)<0,  
\end{align}
where
\begin{align*}
	R(q)=\dfrac{(q;q^5)_\infty(q^4;q^5)_\infty}{(q^2;q^5)_\infty(q^3;q^5)_\infty}=:\sum_{n=0}^{\infty} d(n)q^n.
\end{align*}

Ramanujan, in his lost notebook \cite[p. 50]{lnb} recorded formulas for $\sum_{n=0}^\infty c(5n+j)q^n$ and $\sum_{n=0}^\infty d(5n+j)q^n$, $0\leq j\leq 4$, which were proved by Andrews \cite{andrews1981ramunujan} (Also see Andrews and Berndt \cite[Chapter 4]{andrewsbcb1}). Andrews used the formulas and a theorem of Gordon \cite{gordon} to give partition-theoretic interpretations of these coefficients, and hence proved that \eqref{thm:rich} and \eqref{thm:rich-2}  hold for all $n$ except $c(2)=c(4)=c(9)=0$, $d(3)=d(8)=0$. Using the quintuple product identity \cite{cooper-qtpi}, Hirschhorn \cite{hirschhorn} found exact $q$-product representations of $\sum_{n=0}^\infty c(5n+j)q^n$ and $\sum_{n=0}^\infty d(5n+j)q^n$, $0\leq j\leq 4$, and concluded the periodicity of the signs of the coefficients $c(n)$ and $d(n)$, with two more exceptions, namely $d(13)=d(23)=0$. 

There are results on periodicity of the signs of the coefficients of certain infinite products in more general settings. For $1\leq r,s<m$, define
\begin{align*}
	F_{m,r,s}(q):= \frac{(q^r;q^m)_{\infty}(q^{m-r};q^m)_{\infty}}{(q^s;q^m)_{\infty}(q^{m-s};q^m)_{\infty}}.
\end{align*}
Note that $F_{5,1,2}(q)=R(q)$. In 1988, Ramanathan \cite{ramanathan} proved the following result.
\begin{theorem}\label{thm:ramanathan}
	Suppose $\gcd(m,r)=1$. Let
	\begin{align*}
		F_{m,2r,r}(q)=\sum_{n=0}^{\infty}b(n)q^n.
	\end{align*}
	If $\gcd(m, 6) = 1$, the signs of the $b(n)$'s are periodic with period $m$.
\end{theorem}
Chan and Yesilyurt \cite{chanyesilhurt} proved Theorem \ref{thm:ramanathan} using the ideas of Hirschhorn \cite{hirschhorn} in a more general setting and without the condition $\gcd(m,6)=1$. They also studied the periodicity of infinite products for some other values of $r$ and $s$. Chern and Tang \cite{cherntang} studied the sign patterns of certain $q$-products related to the Rogers-Ramanujan continued fraction, namely, $R(q)R^2(q^2)$, $R^2(q)/R(q^2)$ and their reciprocals. For work on the periodicity of coefficients of similar infinite products, one can look at the following non-exhaustive list of papers: Dou and Xiao \cite{douxiao}, Hirschhorn \cite{hirch2001}, Lin \cite{lin}, Tang \cite{tang}, Tang and Xia \cite{tangxia}, Xia and Yao \cite{xiayao}, and Xia and Zhou \cite{xiazhou}.

In his second notebook \cite[p. 289]{nb} and the lost notebook \cite[p. 365]{lnb}, Ramanujan recorded the following elegant identity: 
\begin{align}	R^5(q)=R(q^5)\cdot\frac{1-2qR(q^5)+4q^2R^2(q^5)-3q^3R^3(q^5)+q^4R^4(q^5)}{1+3qR(q^5)+4q^2R^2(q^5)+2q^3R^3(q^5)+q^4R^4(q^5)}.\label{R_1^5/R_5}
\end{align}
The identity can also be found in his first letter to Hardy written on January 16, 1913. Various proofs of \eqref{R_1^5/R_5} can be found in the literature. For example, see  the papers by Rogers \cite{rogers1921}, Watson \cite{watson}, Ramanathan \cite{ramanathan-acta}, Yi \cite{yi}, and Gugg \cite{gugg-rama}.

In this paper, we investigate the behavior of the signs of the coefficients of the infinite products  $R^5(q)$, $R^5(q)/R(q^5)$, and their reciprocals appearing in \eqref{R_1^5/R_5}. We also find some interesting congruences satisfied by some coefficients. We state our results in the following theorems.
\begin{theorem}\label{periodR} If $A(n)$ is defined by 
	\begin{align*}
		\dfrac{1}{R^5(q)}&=\sum_{n=0}^{\infty}A(n)q^{n},
	\end{align*}then for all nonnegative integers $n$, we have
	\begin{align}
		A(5n+1)>0,\label{R1}\\
		A(5n+2)>0,\label{R2}\\
		A(5n+3)>0,\label{R3}\\
		A(5n+4)<0.\label{R4}
	\end{align}
\end{theorem}

\begin{theorem}\label{periodalpha} If $B(n)$ is defined by 
	\begin{align*}
		R^5(q)&=\sum_{n=0}^{\infty}B(n)q^{n},
	\end{align*}then for all nonnegative integers $n$, we have
	\begin{align}
		B(5n+1)<0,\label{alpha1}\\
		B(5n+2)>0,\label{alpha2}\\
		B(5n+3)<0,\label{alpha3}\\
		B(5n+4)>0.\label{alpha4}
	\end{align}
\end{theorem}

\begin{theorem}\label{periodbeta}
	If $C(n)$ is defined by 
	\begin{align*}
		\dfrac{R^5(q)}{R(q^5)}&=\sum_{n=0}^{\infty}C(n)q^{n},
	\end{align*}then for all nonnegative integers $n$, we have
	\begin{align}
		C(5n)<0, \label{beta0}\\
		C(5n+1)<0, \label{beta1}\\
		C(5n+2)>0, \label{beta2}\\
		C(5n+3)<0, \label{beta3}\\
		C(5n+4)>0,\label{beta4}
	\end{align}
	except $C(0)=1$.
\end{theorem}
\begin{remark} Theorem \ref{periodbeta} shows that the signs of $C(n)$ are periodic with period 5.
\end{remark}

\begin{theorem}\label{periodgamma}If $D(n)$ is defined by 
	\begin{align*}
		\dfrac{R(q^5)}{R^5(q)}&=\sum_{n=0}^{\infty}D(n)q^{n},
	\end{align*}then for all nonnegative integers $n$, we have
	\begin{align}
		D(5n)<0, \label{gamma0}\\
		D(5n+2)>0, \label{gamma2}\\
		D(5n+3)>0, \label{gamma3}\\
		D(5n+4)<0, \label{gamma4}
	\end{align}
	except $D(0)=1$.
\end{theorem}

We state some congruences satisfied by $A(n)$, $B(n)$, $C(n)$, and $D(n)$ in the next theorem.
\begin{theorem}\label{cong:mod3}
	For all $n\geq0$, we have
	\begin{align}
		A(9n+4)&\equiv 0 \pmod3,\label{mod3_1}\\
		B(9n+2)&\equiv 0 \pmod3,\label{mod3_2}\\
		A(16n+13)&\equiv 0 \pmod4,\label{mod4_1}\\
		B(16n+11)&\equiv 0 \pmod4,\label{mod4_2}\\
		A(15n+r)&\equiv0\pmod{15},  \text{where $r\in\{4, 8, 13, 14\}$},\label{mod15_1}\\
		B(15n+r)&\equiv0\pmod{15},  \text{where $r\in\{2, 6, 11, 12\}$},\label{mod15_2}\\
		C(15n+r)&\equiv0\pmod{30},  \text{where $r\in\{3, 13\}$},\label{mod30_1}\\
		D(15n+r)&\equiv0\pmod{30},  \text{where $r\in\{7, 12\}$}.\label{mod30_2}		
	\end{align}
\end{theorem}

We organize the paper as follows. In the next section, we give some preliminary results on Ramanujan's theta functions and $t$-dissections. Note that  a $t$-dissection of a power series $P(q)$ in $q$ is given by $$
P(q) =\sum_{j=0}^{t-1}q^jP_j(q^t),$$
where $P_j$'s are power series in $q$. In Sections \ref{period1}--\ref{congmod3}, we prove Theorems \ref{periodR}--\ref{cong:mod3}. Some concluding remarks and conjectures are presented in the final section of the paper. 
\section{Preliminaries}\label{preliminaries}
To prove our results, we now define Ramanujan's theta functions. Recall that Ramanujan's general theta function $f(a,b)$ \cite[Equation 1.2.1]{Spirit} is defined by 
\begin{align*}
	f(a,b) := \sum_{k=-\infty}^{\infty} a^{\frac{k(k+1)}{2}}b^{\frac{k(k-1)}{2}}, \quad |ab| < 1.
\end{align*}
Three special cases of $f(a,b)$ are:
\begin{align}
	\varphi(q) &:= f(q,q) = (-q;q^2)^2_\infty(q^2;q^2)_\infty=\dfrac{f_2^5}{f_1^2f_4^2},\label{varphi}\\
	\psi(q) &:= f(q,q^3) = \frac{(q^2;q^2)_\infty}{(q;q^2)_\infty}=\dfrac{f_2^2}{f_1},\label{psi}\\
	f(-q) &:= f(-q,-q^2) = (q;q)_\infty=f_1,\notag
\end{align}
where the product representations follow from Jacobi's triple product identity \cite[p. 35, Entry 19]{bcb3},
\begin{align}\label{jtpi}
	f(a,b) = (-a;ab)_\infty(-b;ab)_\infty(ab;ab)_\infty. 
\end{align}
and, for brevity, we set $f_k:=(q^k;q^k)_\infty$. 

\noindent Replacing $q$ by $-q$ in \eqref{varphi} and \eqref{psi}, we have
\begin{align}
	\varphi(-q)&=\frac{f_1^2}{f_2},\label{varphi(-q)}\\
	\psi(-q)&=\frac{f_1f_4}{f_2}.\label{psi(-q)}
\end{align}

In the following lemma, we recall two useful identities.
\begin{lemma}\cite[p. 51 and p. 350]{bcb3}
	We have
	\begin{align}
		\label{f(q,q^5)} f(q,q^5)&=\psi(-q^3)\chi(q),\\
		\label{f(q,q^2)}f(q,q^2)&=\frac{\varphi(-q^3)}{\chi(-q)},
	\end{align}
	where 
	\begin{align}\label{chiq}
		\chi(q):=(-q;q^2)_\infty=\dfrac{f_2^2}{f_1f_4} ~ and ~\chi(-q)=(q;q^2)_\infty=\dfrac{f_1}{f_2}.
	\end{align}
\end{lemma}

In the next three lemmas, we recall some known 2-, 3-, and 5-dissection formulas.
\begin{lemma}\cite[Lemma 2]{matching}\label{lem-2-dis}
	We have
	\begin{align}
		f_1^4&=\dfrac{f_4^{10}}{f_2^2f_8^4}-4q\dfrac{f_2^2f_8^4}{f_4^2}\label{f_1^4},\\
		\frac{1}{f_{1}^2} &= \frac{f_{8}^5}{f_{2}^5 f_{16}^2} + 2 q \frac{f_{4}^2 f_{16}^2 }{f_{2}^5 f_{8}},\label{dis1byf1^2}\\
		\frac{1}{f_{1}^4} &= \frac{f_{4}^{14}}{f_{2}^{14} f_{8}^4} + 4 q \frac{f_{4}^2 f_{8}^4}{f_{2}^{10}},\label{dis1byf1^4}\\
		\frac{f_1}{f_5}&= \frac{f_2f_8f_{20}^3}{f_4f_{10}^3f_{40}}-q \frac{f_4^2f_{40}}{f_8f_{10}^2},\label{disf1byf5}\\
		\frac{f_5}{f_1}&= \frac{f_8f_{20}^2}{f_2^2f_{40}}+q \frac{f_4^3f_{10}f_{40}}{f_2^3f_8f_{20}}\label{disf5byf1}.
	\end{align}
\end{lemma}
\begin{lemma}\cite[p. 49]{bcb3}
	We have
	\begin{align}
		\varphi(q)&=\varphi(q^9)+2qf(q^3,q^{15}),\label{3-dissec-varphi}\\
		\psi(q)&=f(q^3,q^6)+q\psi(q^9).\label{3-dissec-psi}
	\end{align}
\end{lemma}
\begin{lemma}\cite[p. 165]{Spirit}
	We have
	\begin{align}
		f_1&=f_{25}\left(\frac{1}{R(q^5)}-q-q^2R(q^5) \right)\label{f_1}\intertext{and}
		\frac{1}{f_1}&=\frac{f_{25}^5}{f_5^6}\Bigg(\frac{1}{R^4(q^5)}+\frac{q}{R^3(q^5)}+\frac{2q^2}{R^2(q^5)}+\frac{3q^3}{R(q^5)}+5q^4-3q^5R(q^5)\nonumber \\
		&\quad+2q^6R^2(q^5)-q^7R^3(q^5)+q^8R^4(q^5)\Bigg).\label{1-by-f1}
	\end{align}
\end{lemma}

In the next two lemmas, we recall some useful identities involving $R(q)$, $G(q)$, $H(q)$, and $f_k$.
\begin{lemma}\cite[Lemma 7]{gugg-rama}
	We have
	\begin{align}	R^2(q^5)\frac{f_1^2H^5(q)}{f_{25}^2H(q^5)}&=1-2qR(q^5)+4q^2R^2(q^5)-3q^3R^3(q^5)+q^4R^4(q^5),\label{R_1^5/R_5num}\\
		R^2(q^5)\frac{f_1^2G^5(q)}{f_{25}^2G(q^5)}&=1+3qR(q^5)+4q^2R^2(q^5)+2q^3R^3(q^5)+q^4R^4(q^5)\label{R_1^5/R_5den}.
	\end{align}
\end{lemma}
\begin{lemma}\cite[Theorem 7.4.4]{Spirit}, \cite[Eq. (1.22) and Eq. (2.14)]{baruahbegum}, \cite[Eq. (66) and p. 532]{matching}
	We have
	\begin{align}
		\frac{1}{R^5(q)}-q^2R^5(q)&=11q+\frac{f_1^6}{f_5^6}	\label{R1-5}\\
		&=4q+\frac{f_2f_5^5}{f_1f_{10}^5}+8q^2\dfrac{f_1f_{10}^5}{f_2f_5^5}+16q^3\frac{f_1^2f_{10}^{10}}{f_2^2f_5^{10}},\label{baruah-begum}\\
		\frac{1}{R^3(q)R(q^2)}+q^2R^3(q)R(q^2)&=2q+\frac{f_2f_5^5}{f_1f_{10}^5}+4q^2\dfrac{f_1f_{10}^5}{f_2f_5^5},\label{x^3y}\\	
		\frac{R(q^2)}{R^2(q)}+\frac{R^2(q)}{R(q^2)}&=2\frac{f_2^2f_{10}^{10}}{f_4f_5^8f_{20}^3}+8q^2\frac{f_1f_4f_{10}f_{20}^3}{f_2f_5^5},\label{x^2/y}\\
		\frac{1}{R(q)R^2(q^2)}+q^2R(q)R^2(q^2)&=\frac{f_2^3f_{10}^5}{f_1f_4f_5^3f_{20}^3}+4q^2\frac{f_4f_{20}^3}{f_{10}^4}.\label{xy^2+q^2/xy^2}
	\end{align}
\end{lemma}
\section{Proof of Theorem \ref{periodR}}\label{period1}
\noindent \emph{Proof of \eqref{R1}}. Setting
\begin{align}
	N:=1-2qR(q^5)+4q^2R^2(q^5)-3q^3R^3(q^5)+q^4R^4(q^5)\label{N}\intertext{and}D:=1+3qR(q^5)+4q^2R^2(q^5)+2q^3R^3(q^5)+q^4R^4(q^5),\label{D}
\end{align}
we rewrite \eqref{R_1^5/R_5}, \eqref{R_1^5/R_5num}, and \eqref{R_1^5/R_5den} as
\begin{align}
	\frac{R^5(q)}{R(q^5)}&=\frac{N}{D},\label{N/D}\\
	N&=R^2(q^5)\frac{f_1^2H^5(q)}{f_{25}^2H(q^5)},\label{NH}\\
	D&=R^2(q^5)\frac{f_1^2G^5(q)}{f_{25}^2G(q^5)}.\label{NG}
\end{align}
Multiplying \eqref{NH} and \eqref{NG}, we have
\begin{align}
	ND&=R^4(q^5)\frac{f_1^4}{f_{25}^4}\frac{G^5(q)H^5(q)}{G(q^5)H(q^5)}.\label{ND1}
\end{align}

Now, 
\begin{align}
	H(q)G(q)&=\dfrac{1}{(q;q^5)_\infty (q^2;q^5)_\infty(q^3;q^5)_\infty(q^4;q^5)_\infty}=\dfrac{f_5}{f_1}.\label{GH}
\end{align}
Employing \eqref{GH} in \eqref{ND1}, we have
\begin{align}
	ND&=\frac{f_5^6R^4(q^5)}{f_1f_{25}^5}.\label{ND}
\end{align}

From \eqref{N/D} and \eqref{ND}, we have
\begin{align*}
	\frac{1}{R^5(q)}=\frac{D^2}{R(q^5)ND}=\dfrac{f_1f_{25}^5D^2}{f_5^6R^5(q^5)}.
\end{align*}
By \eqref{f_1} and \eqref{D}, we rewrite the above identity as
\begin{align}
	\frac{1}{R^5(q)}=\sum_{n=0}^{\infty}A(n)q^n&=\frac{f_{25}^6}{f_5^6R^5(q^5)}(1+3qR(q^5)+4q^2R^2(q^5)+2q^3R^3(q^5)+q^4R^4(q^5))^2\nonumber\\&\quad\times\left(\frac{1}{R(q^5)}-q-q^2R(q^5)\right).\label{R(n)}
\end{align}
Extracting the terms of the form $q^{5n+1}$ from \eqref{R(n)}, dividing both sides by $q$, and then replacing $q^5$ by $q$, we have
\begin{align*}
	\sum_{n=0}^{\infty}A(5n+1)q^n=\frac{f_5^6}{f_1^6}\left(\frac{5}{R^5(q)}-40q\right).
\end{align*}
Employing \eqref{R1-5} in the above, we have
\begin{align*}
	\sum_{n=0}^{\infty}A(5n+1)q^n&=5+15q\frac{f_5^6}{f_1^6}+5q^2\frac{f_5^6R^5(q)}{f_1^6}\\
	&=5+15q\frac{f_5^6}{f_1^6}+\frac{5q^2}{(q;q^5)_{\infty}(q^2;q^5)_{\infty}^{11}(q^3;q^5)_{\infty}^{11}(q^4;q^5)_{\infty}},
\end{align*}
from which it easily follows that, for $n\geq0$, $A(5n+1)>0$. This proves \eqref{R1}.

\noindent\emph{Proof of \eqref{R2}}. Extracting the terms of the form $q^{5n+2}$ from \eqref{R(n)}, dividing both sides by $q^2$, and then replacing $q^5$ by $q$, we have
\begin{align}\label{match1}
	\sum_{n=0}^{\infty}A(5n+2)q^n&=10\frac{f_5^6}{f_1^6}\left(\frac{1}{R^4(q)}-3qR(q)\right).
\end{align}

Now, let $b_{25}(n)$ denote the number of 25-regular partitions of $n$, that is, the number of partitions in which parts are not divisible by 25. The generating function of $b_{25}(n)$ is given by
\begin{align*}
	\sum_{n=0}^\infty b_{25}(n)q^n=\frac{f_{25}}{f_1}.
\end{align*}
Employing \eqref{1-by-f1} in the above and then extracting the terms involving $q^{5n}$, we find that
\begin{align}\label{b-25-5n}
	\sum_{n=0}^{\infty}b_{25}(5n)q^n=\frac{f_5^6}{f_1^6}\left(\frac{1}{R^4(q)}-3qR(q)\right).
\end{align}
It follows from \eqref{match1} and \eqref{b-25-5n} that
\begin{align*}
	\sum_{n=0}^{\infty}A(5n+2)q^n&=10\sum_{n=0}^\infty b_{25}(5n)q^n.
\end{align*}
As $b_{25}(5n)>0$ for all $n\geq0$, it readily follows from the above that  $A(5n+2)>0$ for all $n\geq0$, which is \eqref{R2}.

\noindent\emph{Proof of \eqref{R3}}. Extracting the terms of the form $q^{5n+3}$ from \eqref{R(n)}, dividing both sides by $q^3$, and then replacing $q^5$ by $q$, we have
\begin{align}
	\sum_{n=0}^{\infty}A(5n+3)q^n
	&=5\frac{f_5^6}{f_1^6}\left(\frac{1}{R^3(q)}-3qR^2(q)\right).\label{r(5n+3)}
\end{align}

Now, define $F(n)$ by
\begin{align}
	\sum_{n=0}^{\infty}F(n)q^n&=\frac{f_{25}R(q^5)}{f_1}\label{Fn}\\
	&=\frac{1}{(q^{1,2,3,4,6,7,8,9,10,10,11,12,13,14,15,15,16,17,18,19,21,22,23,24};
		q^{25})_\infty},\notag
\end{align}
where $(q^{a_1,a_2,\ldots,a_j};q^k)_\infty:=(q^{a_1};q^k)_\infty (q^{a_2};q^k)_\infty\cdots (q^{a_j};q^k)_\infty.$
Combinatorially, $F(n)$ counts the number of partitions of $n$ into parts not congruent to 0 or $\pm5$ modulo 25 and parts congruent to $\pm10$ modulo 25 have two colors. Clearly, $F(n)>0$ for all $n\geq0$.

Now, employing \eqref{1-by-f1} in \eqref{Fn} and then extracting the terms involving $q^{5n}$, we find that
\begin{align}
	\sum_{n=0}^{\infty}F(5n)q^n=\frac{f_5^6}{f_1^6}\left(\frac{1}{R^3(q)}-3qR^2(q)\right).\label{Fn1}
\end{align}
From \eqref{r(5n+3)} and \eqref{Fn1}, we have
\begin{align*}
	\sum_{n=0}^{\infty}A(5n+3)q^n
	&=5\sum_{n=0}^{\infty}F(5n)q^n,
\end{align*}
from which it follows that $A(5n+3)=5F(5n)$. Now, the positivity of $F(n)$ implies that $A(5n+3)>0$ for all $n\geq0$, which is \eqref{R3}. 

\noindent\emph{Proof of \eqref{R4}}. Extracting the terms of the form $q^{5n+4}$ from \eqref{R(n)}, dividing both sides by $q^4$, and then replacing $q^5$ by $q$, we find that
\begin{align}
	&\sum_{n=0}^{\infty}A(5n+4)q^n\notag\\
	&=-5\frac{f_5^6}{f_1^6}\left(\frac{3}{R^2(q)}+qR^3(q)\right)\label{r(5n+4)}\\
	&=-5\left(\frac{3}{(q,q^4;q^5)_{\infty}^8(q^2,q^3;q^5)_{\infty}^4}+\frac{q}{(q,q^4;q^5)_{\infty}^3(q^2,q^3;q^5)_{\infty}^9}\right),\nonumber
\end{align}
which clearly implies that $A(5n+4)<0$ for all $n\geq0$. 

\section{Proof of Theorem \ref{periodalpha}}\label{period2}
\noindent\emph{Proof of \eqref{alpha1}}. From \eqref{N/D} and \eqref{ND}, we have
\begin{align*}
	R^5(q)=\frac{R(q^5)N^2}{ND}=\frac{f_1f_{25}^5N^2}{f_5^6R^3(q^5)}.
\end{align*}
With the aid of \eqref{f_1}, \eqref{N} and \eqref{ND}, we rewrite the above as
\begin{align}
	R^5(q)=\sum_{n=0}^{\infty}B(n)q^n&=\frac{f_{25}^6}{f_5^6R^3(q^5)}\left(1-2qR(q^5)+4q^2R^2(q^5)-3q^3R^3(q^5)+q^4R^4(q^5)\right)^2\nonumber\\&\quad\times\left(\frac{1}{R(q^5)}-q-q^2R(q^5)\right).\label{alpha(n)}
\end{align}
Extracting the terms of the form $q^{5n+1}$ from both sides of the above, dividing both sides by $q$, and then replacing $q^5$ by $q$, we find that
\begin{align}
	\sum_{n=0}^{\infty}B(5n+1)q^n&=-5\frac{f_5^6}{f_1^6}\left(\frac{1}{R^3(q)}-3qR^2(q)\right).\label{alpha5n+1}
\end{align}
From \eqref{r(5n+3)} and \eqref{alpha5n+1}, we have
\begin{align}
	\sum_{n=0}^{\infty}B(5n+1)q^n&=-\sum_{n=0}^{\infty}A(5n+3)q^n,
\end{align}
which implies that for all $n$, $B(5n+1)=-A(5n+3)$. Therefore, \eqref{R3} implies \eqref{alpha1}.

\noindent\emph{Proof of \eqref{alpha2}}. Extracting the terms of the form $q^{5n+2}$ from both sides of \eqref{alpha(n)}, dividing both sides by $q^2$, and then replacing $q^5$ by $q$, we obtain
\begin{align}
	\sum_{n=0}^{\infty}B(5n+2)q^n&=5\frac{f_5^6}{f_1^6}\left(\frac{3}{R^2(q)}+qR^3(q)\right),\label{alpha5n+2}
\end{align}
which, by \eqref{r(5n+4)}, implies that
\begin{align}
	\sum_{n=0}^{\infty}B(5n+2)q^n&=-\sum_{n=0}^{\infty}A(5n+4)q^n.
\end{align}
It follows that $B(5n+2)=-A(5n+4)$. Thus, \eqref{R4} implies \eqref{alpha2}.

\noindent\emph{Proof of \eqref{alpha3}}. Extracting the terms of the form $q^{5n+3}$ from both sides of \eqref{alpha(n)}, dividing both sides by $q^3$, and then replacing $q^5$ by $q$, we find that
\begin{align}
	\sum_{n=0}^{\infty}B(5n+3)q^n&=-10\frac{f_5^6}{f_1^6}\left(\frac{3}{R(q)}+qR^4(q)\right)\label{B5n3}\\
	&=-10\left(\frac{3}{(q,q^4;q^5)_\infty^7(q^2,q^3;q^5)_\infty^5}+\frac{q}{(q,q^4;q^5)_\infty^2(q^2,q^3;q^5)_\infty^{10}}\right),\nonumber
\end{align}
which readily implies that $B(5n+3)<0$ for all $n\geq0$, which is \eqref{alpha3}.

\noindent\emph{Proof of \eqref{alpha4}}. Extracting the terms of the form $q^{5n+4}$ from both sides of \eqref{alpha(n)}, dividing both sides by $q^4$, and then replacing $q^5$ by $q$, we obtain
\begin{align*}
	\sum_{n=0}^{\infty}B(5n+4)q^n&=\frac{f_5^6}{f_1^6}\left(40+5qR^5(q)\right)\\
	&=40\frac{f_5^6}{f_1^6}+\frac{5q}{(q,q^4;q^5)_\infty(q^2,q^3;q^5)_\infty^{11}},
\end{align*}
which implies that $B(5n+4)>0$ for all $n\geq0$.

\section{Proof of Theorem \ref{periodbeta}}\label{period3}
\noindent\emph{Proof of \eqref{beta0}}. From \eqref{N/D}, we have
\begin{align*}
	\frac{R^5(q)}{R(q^5)}=\frac{N^2}{ND}=\frac{f_1f_{25}^5N^2}{f_5^6R^4(q^5)}.
\end{align*}
Invoking \eqref{f_1}, \eqref{N} and \eqref{ND} in the above, we obtain
\begin{align}
	\frac{R^5(q)}{R(q^5)}=\sum_{n=0}^{\infty}C(n)q^n&=\frac{f_{25}^6}{f_5^6R^4(q^5)}(1-2qR(q^5)+4q^2R^2(q^5)-3q^3R^3(q^5)+q^4R^4(q^5))^2\nonumber\\&\quad\times\left(\frac{1}{R(q^5)}-q-q^2R(q^5)\right).\label{beta(n)}
\end{align}
Extracting the terms of the form $q^{5n}$ and then replacing $q^5$ by $q$, we have
\begin{align}\label{C5n}
	\sum_{n=0}^{\infty}C(5n)q^n=\frac{f_5^6}{f_1^6}\left(\frac{1}{R^5(q)}-36q-q^2R^5(q)\right).
\end{align}
With the aid of \eqref{R1-5}, the last equation can be recast as
\begin{align*}
	\sum_{n=0}^{\infty}C(5n)q^n=1-25q\frac{f_5^6}{f_1^6}
\end{align*}
which clearly implies $C(5n)<0$ for all $n\geq1$, which is \eqref{beta0}.

\noindent\emph{Proof of \eqref{beta1}}. Extracting the terms of the form $q^{5n+1}$ from \eqref{beta(n)}, dividing both sides by $q$ and then replacing $q^5$ by $q$, we have
\begin{align}
	\sum_{n=0}^{\infty}C(5n+1)q^n&=-5\frac{f_5^6}{f_1^6}\left(\frac{1}{R^4(q)}-3qR(q)\right).\label{match2}
\end{align}
It follows from \eqref{match1} and \eqref{match2} that, for all $n\geq0$, we have $2C(5n+1)=-A(5n+2)$. Thus, \eqref{beta1} follows by \eqref{R2}.

\noindent\emph{Proof of \eqref{beta2}}. Extracting the terms of the form $q^{5n+2}$ from \eqref{beta(n)}, dividing both sides by $q^2$ and then replacing $q^5$ by $q$, we have
\begin{align*}
	\sum_{n=0}^{\infty}C(5n+2)q^n&=\frac{f_5^6}{f_1^6}\left(\frac{15}{R^3(q)}+5qR^2(q)\right)\\
	&=\frac{15}{(q,q^4;q^5)_\infty^9(q^2,q^3;q^5)_\infty^3}+\frac{5q}{(q,q^4;q^5)_\infty^4(q^2,q^3;q^5)_\infty^8},
\end{align*}
which readily yields $C(5n+2)>0$ for $n\geq0$, which is \eqref{beta2}.

\noindent\emph{Proof of \eqref{beta3}}. Extracting the terms of the form $q^{5n+3}$ from \eqref{beta(n)}, dividing both sides by $q^3$ and then replacing $q^5$ by $q$, we have
\begin{align}
	\sum_{n=0}^{\infty}C(5n+3)q^n&=-10\frac{f_5^6}{f_1^6}\left(\frac{3}{R^2(q)}+qR^3(q)\right).\label{beta5n+3}
\end{align}
From \eqref{r(5n+4)} and \eqref{beta5n+3}, we have $C(5n+3)=2A(5n+4)$ for all $n\geq0$. Thus, \eqref{R4} implies \eqref{beta3}.

\noindent\emph{Proof of \eqref{beta4}}. Extracting the terms of the form $q^{5n+4}$ from \eqref{beta(n)}, dividing both sides by $q^4$, and then replacing $q^5$ by $q$, we have
\begin{align*}
	\sum_{n=0}^{\infty}C(5n+4)q^n&=\frac{f_5^6}{f_1^6}\left(\frac{40}{R(q)}+5qR^4(q)\right)\\
	&=\frac{40}{(q,q^4;q^5)_\infty^7(q^2,q^3;q^5)_\infty^5}+\frac{5q}{(q,q^4;q^5)_\infty^2(q^2,q^3;q^5)_\infty^{10}},
\end{align*}
which yields $C(5n+4)>0$ for all $n\geq0$, which is \eqref{beta4}.
\section{Proof of Theorem \ref{periodgamma}}\label{period4}
\noindent\emph{Proof of \eqref{gamma0}}. From \eqref{N/D}, we have
\begin{align*}
	\frac{R(q^5)}{R^5(q)}=\frac{D^2}{ND}=\frac{f_1f_{25}^5D^2}{f_5^6R^4(q^5)}.
\end{align*}
Invoking \eqref{f_1}, \eqref{D} and \eqref{ND} in the last equation, we arrive at
\begin{align}
	\frac{R(q^5)}{R^5(q)}=\sum_{n=0}^{\infty}D(n)q^n&=\frac{f_{25}^6}{f_5^6R^4(q^5)}(1+3qR(q^5)+4q^2R^2(q^5)+2q^3R^3(q^5)+q^4R^4(q^5))^2\nonumber\\&\quad\times\left(\frac{1}{R(q^5)}-q-q^2R(q^5)\right).\label{gamma(n)}
\end{align}
Extracting the terms of the form $q^{5n}$ and then replacing $q^5$ by $q$, we have
\begin{align*}
	\sum_{n=0}^{\infty}D(5n)q^n=\frac{f_5^6}{f_1^6}\left(\frac{1}{R^5(q)}-36q-q^2R^5(q)\right).
\end{align*}
Comparing the above identity with \eqref{C5n}, we see that $D(5n)=C(5n)$ for all $n\geq0$. Thus,  \eqref{beta0} implies \eqref{gamma0}.

\noindent\emph{Proof of \eqref{gamma2}}. Extracting the terms of the form $q^{5n+2}$ from \eqref{gamma(n)}, dividing both sides by $q^2$ and then replacing $q^5$ by $q$, we have
\begin{align}
	\sum_{n=0}^{\infty}D(5n+2)q^n&=10\frac{f_5^6}{f_1^6}\left(\frac{1}{R^3(q)}-3qR^2(q)\right).\label{gamma5n+2}
\end{align}
It follows from the above identity and \eqref{r(5n+3)} that $D(5n+2)=2A(5n+3)$ for all $n\geq0$. Therefore, \eqref{gamma2} follows by \eqref{R3}.

\noindent\emph{Proof of \eqref{gamma3}}.
Extracting the terms of the form $q^{5n+3}$ from \eqref{gamma(n)}, dividing both sides by $q^3$ and then replacing $q^5$ by $q$, we have
\begin{align}
	\sum_{n=0}^{\infty}D(5n+3)q^n&=\frac{f_5^6}{f_1^6}\left(\frac{5}{R^2(q)}-15qR^3(q)\right)\nonumber\\
	&=\frac{5}{(q,q^4;q^5)^3_{\infty}(q^2,q^3;q^5)^3_{\infty}}\left(\frac{1}{(q,q^4;q^5)^5_{\infty}}-\frac{3q}{(q^2,q^3;q^5)^5_{\infty}}\right).\label{gamma31}
\end{align}
In order to prove \eqref{gamma3}, that is, $D(5n+3)>0$ for all nonnegtive integers $n$, it is enough to show that coefficients in the expansion of 
$$\left(\frac{1}{(q,q^4;q^5)^5_{\infty}}-\frac{3q}{(q^2,q^3;q^5)^5_{\infty}}\right)$$are positive.

To that end, define the sequences $(\alpha(n))$ and $(\beta(n))$ by
\begin{align*}
	\sum_{n=0}^{\infty}\alpha(n)q^n:=\frac{1}{(q,q^4;q^5)_\infty^5}
	\intertext{and}
	\sum_{n=0}^{\infty}\beta(n)q^n:=\frac{1}{(q^2,q^3;q^5)_\infty^5}.
\end{align*}
Combinatorially, $\alpha(n)$ counts the number of partitions 5-tuples of $n$ with parts congruent to $\pm1\pmod5$ and $\beta(n)$ counts the number of partitions 5-tuples of $n$ with parts congruent to $\pm2\pmod5$. We first show that
\begin{align}
	\alpha(n)>3\beta(n-1) ~\text{for $n\geq3$}.\label{gamma32}
\end{align}
We will follow a procedure similar to Tang \cite[pp. 447--448]{tang}. Define the following sets:
\begin{align*}
	&R_1(n):=\{\mathcal{P}|\text{$\mathcal{P}$ is a partition of $n$ and all parts are $\equiv\pm 1\pmod5$}\},\\
	&R_2(n):=\{\mathcal{P}|\text{$\mathcal{P}$ is a partition of $n$ and all parts are $\equiv\pm 2\pmod5$}\},\\
	&\mathcal{E}_n:=\{\pi=(\pi_1,\pi_2,\pi_3,\pi_4,\pi_5)|\text{all parts in partitions $\pi_i, 1\leq i\leq5$, are $\equiv\pm 2 \pmod5$,}\\&\quad\quad\quad \sum_{i=1}^{5}s(\pi_i)=n\},\\
	&\mathcal{S}_n:=\{\pi=(\pi_1,\pi_2,\pi_3,\pi_4,\pi_5)|\text{all parts in partition $\pi_i, 1\leq i\leq5,$ are $\equiv\pm 1 \pmod5$,}\\&\quad\quad\quad \sum_{i=1}^{5}s(\pi_i)=n\},
\end{align*}
where $s(\pi_i)$ denotes the sum of all parts in partition $\pi_i$.

Define the map $\tau: R_2(n)\to R_1(n)$ by $\tau(\pi)=\lambda$, where $\lambda$ is a partition obtained by subtracting 1 from each of the parts congruent to $2\pmod5$ and adding 1 to each of the the parts congruent to 3$\pmod5$. In case the number of parts congruent to $2\pmod5$ and $3\pmod5$ are unequal, then we do the following for the extra partitions:
\begin{enumerate}
	\item If the number of parts congruent to $2\pmod5$ are more than the number of parts congruent to $3\pmod5$, then we write each remaining part as a part congruent to $1\pmod5$ + 1.
	\item If the number of parts congruent to $3\pmod5$ are more than the number of parts congruent to $2\pmod5$, then we write each remaining part as a part congruent to $1\pmod5$ + 1 + 1.
\end{enumerate}
Clearly, $\lambda\in R_1(n)$. Furthermore, $\pi_1\neq\pi_2\iff\tau(\pi_1)\neq\tau(\pi_2)$.

Again, define $\tilde{\tau}:\mathcal{E}_n\to\mathcal{S}_n$ by
\begin{align*}
	\tilde{\tau}=(\tau(\pi_1),\tau(\pi_2),\tau(\pi_3),\tau(\pi_4),\tau(\pi_5))=\lambda^\prime.
\end{align*} 
If $\pi\in\mathcal{E}_{n-1}$, then $\lambda^\prime\in\mathcal{S}_{n-1}$. We now add a part of size one to any one of the components of $\lambda^\prime$ and let $\tilde{\lambda^\prime}$ be the new partition. Clearly, $\tilde{\lambda^\prime}\in\mathcal{S}_n$. Since, for the partition 5-tuple, there are 5 choices to append this part of size one, so $\alpha(n)\geq5\beta(n-1)$. Furthermore, since $\beta(n)\neq0$ for all $n\geq2$, we see that $\alpha(n)\geq5\beta(n-1)>3\beta(n-1)$ for $n\geq3$. Thus, \eqref{gamma32} holds. Using \eqref{gamma32} and the easily checked facts $D(3)=5$, $D(8)=25$, and $D(13)=155$, in \eqref{gamma31}, we readily arrive at \eqref{gamma3}.

\noindent\emph{Proof of \eqref{gamma4}}. Extracting the terms of the form $q^{5n+4}$ from \eqref{gamma(n)}, dividing both sides by $q^4$ and then replacing $q^5$ by $q$, we have
\begin{align*}
	\sum_{n=0}^{\infty}D(5n+4)q^n&=-5\frac{f_5^6}{f_1^6}\left(\frac{3}{R(q)}+qR^4(q)\right).
\end{align*}
From the above identity and \eqref{B5n3}, we see that $2D(5n+4)=B(5n+3)$ for all $n\geq0$. Therefore, \eqref{alpha3} implies \eqref{gamma4}.
\section{Proof of Theorem \ref{cong:mod3}}\label{congmod3}
\noindent\emph{Proofs of \eqref{mod3_1} and \eqref{mod3_2}}. From the definitions of $A(n)$ and $B(n)$ given in Theorem \ref{periodR} and Theorem \ref{periodalpha} respectively, the identity \eqref{R1-5} may be rewritten as
\begin{align}\label{mod3-a}
	\sum_{n=0}^{\infty}A(n)q^{n}-\sum_{n=0}^{\infty}B(n)q^{n+2}&=11q+\frac{f_1^6}{f_5^6}.
\end{align}

Now, with the aid of the binomial theorem, it can be easily shown that 
\begin{align}\label{cong-mod3}
	f_1^3 \equiv f_3 \pmod3.
\end{align}
Therefore, from \eqref{mod3-a}, we have
\begin{align}
	\sum_{n=0}^{\infty}A(n)q^{n}-\sum_{n=0}^{\infty}B(n)q^{n+2}&\equiv 2q+\frac{f_3^2}{f_{15}^2}\pmod3.\label{R(n)-alpha(n)}
\end{align}

Furthermore, by \eqref{GH}, we have
\begin{align*}
	\sum_{n=0}^{\infty}A(n)q^{n}+\sum_{n=0}^{\infty}B(n)q^{n+2}&=\dfrac{1}{R^5(q)}+q^2R^5(q)\\
	&=\frac{G^5(q)}{H^5(q)}+q^2\frac{H^5(q)}{G^5(q)}\\
	&=\frac{f_1^5}{f_5^5}\left(G^{10}(q)+q^2H^{10}(q)\right),
\end{align*}
which implies that
\begin{align}\label{mod3-1-i}
	\sum_{n=0}^{\infty}A(n)q^{n}+\sum_{n=0}^{\infty}B(n)q^{n+2}
	&\equiv\frac{f_1^2f_3}{f_5^2f_{15}}\left(G(q)G^3(q^3)+q^2H(q)H^3(q^3)\right)\pmod3,
\end{align}
where we also applied the congruences
\begin{align*}
	G^3(q) \equiv G(q^3) \pmod3 ~\textup{and}~H^3(q) \equiv H(q^3) \pmod3.
\end{align*}

Now, recall from \cite[Theorem 4.1(ii)]{gugg-jnt} that
\begin{align*}
	G(q)G^3(q^3)+q^2H(q)H^3(q^3)=\frac{f_3^2}{f_1f_9}.
\end{align*}
Using the above identity in \eqref{mod3-1-i} and then employing \eqref{cong-mod3}, we find that 
\begin{align}
	\sum_{n=0}^{\infty}A(n)q^{n}+\sum_{n=0}^{\infty}B(n)q^{n+2}&\equiv\frac{f_1f_5}{f_{15}^2}\pmod3.\label{R(n)+alpha(n)}
\end{align}
Adding \eqref{R(n)-alpha(n)} and \eqref{R(n)+alpha(n)}, we have
\begin{align}
	2\sum_{n=0}^{\infty}A(n)q^{n}\equiv 2 q+\frac{f_3^2}{f_{15}^2}+\frac{f_1f_5}{f_{15}^2}\pmod3.\label{2R(n)mod3}
\end{align}

Again, from \cite[p. 509]{ndb-boruah}, we recall that
\begin{align}
	f_1f_5&=\varphi(q^5)\psi(q^2)-q\varphi(q)\psi(q^{10}).\label{eq:ndbboruah}
\end{align}
Employing \eqref{3-dissec-varphi} and \eqref{3-dissec-psi} in the above, we have
\begin{align}
	f_1f_5&=\left(\varphi(q^{45})+2q^5f(q^{15},q^{75})\right)\left(f(q^6,q^{12})+q^2\psi(q^{18})\right)-q\left(\varphi(q^9)+2qf(q^9,q^{15})\right)\nonumber\\&\quad\times\left(f(q^{30},q^{60})+q^{10}\psi(q^{90})\right).\label{f_11f_5}
\end{align}
Invoking \eqref{f_11f_5} in \eqref{2R(n)mod3}, and then extracting the terms involving $q^{3n+1}$, we find that
\begin{align*}
	2\sum_{n=0}^{\infty}A(3n+1)q^{n}&\equiv2+\frac{2q^2f(q^5,q^{25})\psi(q^6)-\varphi(q^3)f(q^{10},q^{20})}{f_5^2}\pmod3,
\end{align*}
which, by \eqref{f(q,q^5)} and \eqref{f(q,q^2)}, can be rewritten as
\begin{align*}
	2\sum_{n=0}^{\infty}A(3n+1)q^{n}&\equiv2+2q^2\frac{\psi(q^6)\psi(-q^{15})\chi(q^5)}{f_5^2}-\frac{\varphi(q^3)\varphi(-q^{30})}{f_5^2\chi(-q^{10})}\pmod3.
\end{align*}
Employing \eqref{chiq}, \eqref{varphi(-q)}, and \eqref{psi(-q)}, in the above, we find that
\begin{align*}
	2\sum_{n=0}^{\infty}A(3n+1)q^{n}&\equiv2+2q^2\frac{\psi(q^6)\psi(-q^{15})f_{10}^2}{f_{15}f_{20}}-\frac{\varphi(q^3)\varphi(-q^{30})f_{20}f_5}{f_{15}f_{10}}\\
	&\equiv2+2q^2\frac{\psi(q^6)\psi(-q^{15})}{f_{15}}\varphi(-q^{10})-\frac{\varphi(q^3)\varphi(-q^{30})}{f_{15}}\psi(-q^{5})\pmod3.
\end{align*}
With the aid of \eqref{3-dissec-varphi} and \eqref{3-dissec-psi}, the above can be written as
\begin{align*}
	2\sum_{n=0}^{\infty}A(3n+1)q^{n}&\equiv2+2q^2\frac{\psi(q^6)\psi(-q^{15})}{f_{15}}\left(\varphi(-q^{90})-2q^{10}f(-q^{30},-q^{150})\right)\\&\quad-\frac{\phi(q^3)\psi(-q^{30})}{f_{15}}\left(f(-q^{15},q^{30})-q^{5}\psi(-q^{45})\right)\pmod3.
\end{align*}
Comparing the coefficients of $q^{3n+1}$ for $n\geq0$, we obtain
\begin{align}
	A(9n+4)\equiv0\pmod3,\label{R(9n+4)cong0mod3}
\end{align}
which is \eqref{mod3_1}.

Again, extracting the terms involving $q^{3n+4}$ from both sides of \eqref{R(n)-alpha(n)}, we find that
\begin{align*}
	A(3n+4)\equiv B(3n+2)\pmod3.
\end{align*}
Replacing $n$ by $3n$ in the above, and then invoking \eqref{R(9n+4)cong0mod3}, we have
\begin{align*}
	A(9n+4)\equiv B(9n+2)\equiv0\pmod3,
\end{align*}
which is \eqref{mod3_2}.

\noindent\emph{Proof of \eqref{mod4_1}}. At first, from \eqref{baruah-begum} and the definitions of $A(n)$ and $B(n)$ given in Theorem \ref{periodR} and Theorem \ref{periodalpha}, we have
\begin{align}
	\sum_{n=0}^{\infty}A(n)q^n-\sum_{n=0}^{\infty}B(n)q^{n+2}
	&\equiv\frac{f_2f_5^5}{f_1f_{10}^5}+4q\pmod{8}.\label{x^5-q^2/x^5}
\end{align}

Next,
\begin{align*}
	\sum_{n=0}^{\infty}A(n)q^n+\sum_{n=0}^{\infty}B(n)q^{n+2}&=\frac{1}{R^5(q)}+q^2R^5(q)\\
	&=-\frac{1}{R(q)R^2(q^2)}-q^2R(q)R^2(q^2)+\left(\frac{R^2(q)}{R(q^2)}+\frac{R(q^2)}{R^2(q)}\right)\\&\quad\times\left(\frac{1}{R^3(q)R(q^2)}+q^2R^3(q)R(q^2)\right).\\
\end{align*}
Invoking \eqref{x^3y}--\eqref{xy^2+q^2/xy^2} in the above, we find that
\begin{align}
	\sum_{n=0}^{\infty}A(n)q^n+\sum_{n=0}^{\infty}B(n)q^{n+2}&=-\frac{f_2^3f_{10}^5}{f_1f_4f_5^3f_{20}^3}-4q^2\frac{f_4f_{20}^3}{f_{10}^4}+\left(2\frac{f_2^2f_{10}^{10}}{f_4f_5^8f_{20}^3}+8q^2\frac{f_1f_4f_{10}f_{20}^3}{f_2f_5^5}\right)\nonumber\\&\quad
	\times\left(2q+\frac{f_2f_5^5}{f_1f_{10}^5}+4q^2 \frac{f_1f_{10}^5}{f_2f_5^5}\right)\nonumber\\
	&\equiv\frac{f_2^3f_{10}^5}{f_1f_4f_5^3f_{20}^3}+4q+4q^2f_4f_{20}\pmod{8}.\label{modAn0Bn}
\end{align}
Here and throughout the sequel we use the fact that for integers $k\geq1$ and $\ell\geq1$,
\begin{align*}
	f_k^{2^\ell} \equiv f_{2k}^{2^{\ell-1}} \pmod{2^\ell}.
\end{align*}

We now add \eqref{x^5-q^2/x^5} and \eqref{modAn0Bn}, and then use \eqref{f_1^4}, \eqref{dis1byf1^4}, and \eqref{disf5byf1}. Accordingly, we find that
\begin{align*}
	2\sum_{n=0}^{\infty}A(n)q^n&\equiv\frac{f_2^3f_{10}^5}{f_4f_{20}^3}\cdot \frac{1}{f_5^4}\cdot\frac{f_5}{f_1}+\frac{f_2}{f_{10}^5}\cdot f_5^4\cdot\frac{f_5}{f_1}+4q^2f_4f_{20}\\
	&\equiv\frac{f_2^3f_{10}^5}{f_4f_{20}^3}\left(\frac{f_8f_{20}^2}{f_2^2f_{40}}+q\frac{f_4^3f_{10}f_{40}}{f_2^3f_8f_{20}}\right)\left(\frac{f_{20}^{14}}{f_{10}^{14}f_{40}^4}+4q^5\frac{f_{20}^2f_{40}^4}{f_{10}^{10}}\right)+\frac{f_2}{f_{10}^5}\\&\quad \times\left(\frac{f_8f_{20}^2}{f_2^2f_{40}}+q\frac{f_4^3f_{10}f_{40}}{f_2^3f_8f_{20}}\right)\left(\frac{f_{20}^{10}}{f_{10}^{2}f_{40}^4}-4q^5\frac{f_{10}^2f_{40}^4}{f_{20}^{2}}\right)+4q^2f_4f_{20}\pmod{8}.
\end{align*}
Extracting the terms involving odd powers of $q$, we have
\begin{align*}
	2\sum_{n=0}^{\infty}A(2n+1)q^n&\equiv\frac{f_2^3f_{10}^9}{f_1^2f_4f_5^6f_{20}^3}+\frac{f_2^2f_{10}^{10}}{f_4f_5^8f_{20}^3}-4q^2\frac{f_4f_{20}^3}{f_1f_5^3}+4q^2\frac{f_1f_4f_{10}f_{20}^3}{f_2f_5^5}\\
	&\equiv \frac{f_2^3f_{10}^5}{f_4f_{20}^3}\cdot \left(\frac{f_5}{f_1}\right)^2+\frac{f_2^2f_{10}^{6}}{f_4f_{20}^3}\pmod8.
\end{align*}
Employing \eqref{disf5byf1} in the above and then extracting the terms involving the even powers of $q$, we obtain
\begin{align*}
	2\sum_{n=0}^{\infty}A(4n+1)q^n&\equiv\frac{f_{10}}{f_2}\cdot\left(\frac{f_1}{f_5}\right)^2+\frac{f_4^2f_{10}}{f_2f_{20}^2}\cdot\frac{f_5}{f_1}\cdot f_5^4+q\frac{f_2^5f_{20}^2}{f_4^2f_{10}}\cdot\frac{f_1}{f_5}\cdot\frac{1}{f_1^4}\pmod8.
\end{align*}
Using the identities of Lemma \ref{lem-2-dis} in the above and then extracting  the odd powers of $q$, we find that
\begin{align*}
	2\sum_{n=0}^{\infty}A(8n+5)q^n&\equiv\frac{f_2^5f_{20}}{f_4f_{10}}\cdot\frac{1}{f_1^4}-4q^2 f_4f_{20}^3\cdot\frac{f_1}{f_5}+\frac{f_4f_{10}^5}{f_2f_{20}}\cdot\frac{1}{f_5^4}-4qf_4^3f_{20}\cdot\frac{f_5}{f_1}\\
	&\quad-2f_2f_{10}\pmod8.
\end{align*}
Once again using the identities of Lemma \ref{lem-2-dis} in the above and then extracting  the odd powers of $q$, we have
\begin{align*}
	2\sum_{n=0}^{\infty}A(16n+13)q^n&\equiv4\frac{f_2f_4^4f_{10}}{f_1^5f_5}+4q\frac{f_2^3f_{10}^3f_{20}}{f_4f_5^2}+4q^2\frac{f_2f_{10}f_{20}^4}{f_1f_5^5}-4\frac{f_2^3f_4f_{10}^3}{f_1^2f_{20}}\pmod8
\end{align*}
which implies that
\begin{align*}
	\sum_{n=0}^{\infty}A(16n+13)q^n&\equiv2\frac{f_2f_4^4f_{10}}{f_1^5f_5}+2q\frac{f_2^3f_{10}^3f_{20}}{f_4f_5^2}+2q^2\frac{f_2f_{10}f_{20}^4}{f_1f_5^5}-2\frac{f_2^3f_4f_{10}^3}{f_1^2f_{20}}\pmod4.
\end{align*}
Equivalently,
\begin{align}
	\sum_{n=0}^{\infty}A(16n+13)q^n&\equiv2f_1f_5f_4^3+2qf_2f_{10}^4+2q^2f_1f_5f_{20}^3-2f_8f_{10}\nonumber\\
	&\equiv2f_2f_5(f_2^3-qf_{10}^3)(f_2^3+qf_{10}^3-f_1f_5)\pmod4.\label{A(16n+13final)}
\end{align}

Now, from \eqref{eq:ndbboruah}, we have
\begin{align*}
	f_1f_5&=\varphi(q^5)\psi(q^2)-q\varphi(q)\psi(q^{10})\\
	&=\frac{f_{10}^5}{f_5^2f_{20}^2}\cdot\frac{f_4^2}{f_2}-q\frac{f_{2}^5}{f_{1}^2f_{4}^2}\cdot\frac{f_{20}^2}{f_{10}}\\
	&\equiv f_2^3+q f_{10}^3\pmod2.
\end{align*}
Employing the above congruence in \eqref{A(16n+13final)}, we find that 
\begin{align*}
	\sum_{n=0}^{\infty}A(16n+13)q^n&\equiv2f_2f_5(f_2^3-qf_{10}^3)(f_2^3+qf_{10}^3-(f_2^3+q f_{10}^3))\equiv0\pmod4,
\end{align*}
which is \eqref{mod4_1}.

\noindent\emph{Proof of \eqref{mod4_2}}. From \eqref{baruah-begum}, we have
\begin{align*}
	\sum_{n=0}^{\infty}A(n)q^n-\sum_{n=0}^{\infty}B(n)q^{n+2}&\equiv\frac{f_2f_5}{f_1f_{10}^3}\pmod{4},
\end{align*}
which, by \eqref{disf5byf1}, may be written as
\begin{align*}
	\sum_{n=0}^{\infty}A(n)q^n-\sum_{n=0}^{\infty}B(n)q^{n+2}&\equiv\frac{f_2}{f_{10}^3}\left(\frac{f_8f_{20}^2}{f_2^2f_{40}}+q\frac{f_4^3f_{10}f_{40}}{f_2^3f_8f_{20}}\right)\pmod4.
\end{align*}
Extracting, we have
\begin{align*}
	\sum_{n=0}^{\infty}A(2n+1)q^n-\sum_{n=1}^{\infty}B(2n-1)q^n&\equiv\frac{f_2^3f_{20}}{f_4f_{10}}\cdot\frac{1}{f_1^2}\cdot\frac{1}{f_5^2}\pmod4.
\end{align*}
Employing \eqref{dis1byf1^2} in the above and then extracting, we find that
\begin{align*}
	\sum_{n=0}^{\infty}A(4n+1)q^n-\sum_{n=1}^{\infty}B(4n-1)q^n&\equiv\frac{f_4^5f_{20}^5}{f_2f_{8}^2f_{10}f_{40}^2}\cdot\frac{1}{f_1^2}\cdot\frac{1}{f_5^2}\pmod4.
\end{align*}
Using \eqref{dis1byf1^2} again in the above and then extracting the odd powers of $q$, we obtain
\begin{align*}
	\sum_{n=0}^{\infty}A(8n+5)q^n-\sum_{n=1}^{\infty}B(8n+3)q^n&\equiv2\frac{f_2^7f_8^2 f_{10}^5 f_{20}^3 }{f_1^6 f_4^3 f_5^6 f_{40}^2}+2q^2\frac{f_2^5 f_4^3 f_{10}^7 f_{40}^2 }{f_1^6 f_5^6 f_8^2 f_{20}^3}\\
	&\equiv2\frac{f_2^4f_8^2 f_{10}^2 f_{20}^3 }{f_4^3 f_{40}^2}+2q^2\frac{f_2^2 f_4^3 f_{10}^4 f_{40}^2 }{f_8^2 f_{20}^3}\pmod4,
\end{align*}
which readily implies that
\begin{align*}
	A(16n+13)\equiv B(16n+11)\pmod4.
\end{align*}
Using \eqref{mod4_1}, we arrive at
\begin{align*}
	B(16n+11)\equiv0\pmod{4},
\end{align*}
which is \eqref{mod4_2}.

\noindent\emph{Proofs of \eqref{mod15_1}--\eqref{mod30_2}}.
From \eqref{r(5n+3)}, we have
\begin{align*}
	\sum_{n=0}^{\infty}A(5n+3)q^n&=\frac{f_5^6}{f_1^6}\left(\frac{5}{R^3(q)}-15qR^2(q)\right)\\
	&\equiv5\frac{f_5^6}{f_1^6R^3(q)}\pmod{15}\\
	&\equiv5\frac{f_{15}^2}{f_3^2R(q^3)}\pmod{15}.
\end{align*}
Clearly, the last relation has no terms involving $q^{3n+r}$ for  $r\in\{1, 2\}$. Hence, for all $n\geq0$, we have
\begin{align*} A(15n+8)\equiv A(15n+13)\equiv 0\pmod{15}.\end{align*}

To prove the remaining two congruences in \eqref{mod15_1}, we consider \eqref{r(5n+4)} and proceed exactly as in the above. Thus, we complete the proof of \eqref{mod15_1}.

The proofs of \eqref{mod15_2}--\eqref{mod30_2} are similar to the proof of \eqref{mod15_1}. So we only record the required generating functions for the proofs in the following table.
\begin{center}
	\begin{tabular}{ | m{5cm} | m{5cm}| } 
		\hline
		Congruence & Generating functions\\
		\hline
		\eqref{mod15_2}& \eqref{alpha5n+1} and \eqref{alpha5n+2}\\
		\hline
		\eqref{mod30_1}& \eqref{beta5n+3}\\
		\hline
		\eqref{mod30_2}& \eqref{gamma5n+2}\\
		\hline
	\end{tabular}
\end{center}
\section{Concluding remarks}\label{sec:conclude}
In Theorems \ref{periodR}, \ref{periodalpha}, and \ref{periodgamma}, we have the sign patterns of the coefficients $A(n)$, $B(n)$, and $D(n)$ of $1/R^5(q)$, $R^5(q)$, and $R(q^5)/R^5(q)$, respectively, except $A(5n)$, $B(5n)$, and $D(5n+1)$. Based on numerical observation, we pose the following conjecture.

\begin{conjecture}
	For all integers $n\geq0$, 
	\begin{align*}
		A(5n)&<0,\\
		B (5n)&<0,\\
		D(5n+1)&>0.			
	\end{align*}
\end{conjecture}
An affirmative answer to the above conjecture along with Theorems \ref{periodR}, \ref{periodalpha}, and \ref{periodgamma} will prove that 
the signs of $A(n), B(n)$, and $D(n)$ are periodic with period 5.

There might be more congruences similar to those given in Theorem \ref{cong:mod3}. For example, based on numerical calculations, we propose the following conjecture.
\begin{conjecture} For all integers $n\geq0$,
	\begin{align*} 
		C(27n+18)&\equiv 0 \pmod3,\\
		D(27n+18)&\equiv 0 \pmod3,\\
		C(16n+12)&\equiv 0 \pmod4,\\
		D(16n+12)&\equiv 0 \pmod4,\\
		C(32n+28)&\equiv 0 \pmod8,\\
		D(32n+28)&\equiv 0 \pmod8.
	\end{align*}
\end{conjecture}
\section*{Acknowledgments}

The authors are indebted to the referee for his/her helpful comments. The second author was partially supported by an institutional fellowship for doctoral research from Tezpur University, Assam, India. The author thanks the funding institution.

\end{document}